\documentclass[12pt]{article}

\usepackage{paralist}
\usepackage{amssymb}
\usepackage{amsfonts}
\usepackage{amsmath}
\usepackage{graphicx}
\usepackage{latexsym}
\usepackage{epsfig}

\usepackage[colorlinks=true]{hyperref}
\newtheorem{theorem}{Theorem}

\newtheorem{definition}{Definition}
\newtheorem{example}{Example}

\newtheorem{lemma}{Lemma}

\newtheorem{proposition}{Proposition}

\newtheorem{remark}{Remark}

\title{Zeta-functions of renormalizable sub-Lorenz templates}

\begin{document}

\maketitle
%% Enter the first author's name and address:
\centerline{\scshape Nuno Franco\footnote{Both authors are
partially supported by FCT-Portugal. The first author is partially
supported by FCT grant SFRH/BPD/26354/2006} }
\medskip
{\footnotesize
 %% please put the address of the first author
 \centerline{CIMA-UE and Department of Mathematics, University of \'{E}vora}
   \centerline{Rua Rom\~{a}o Ramalho, 59, 7000-671 \'{E}vora, Portugal, email nmf@uevora.pt}
} %% Do not forget to end the {\footnotesize by the sign }

\medskip

\centerline{\scshape and Lu\'
{i}s Silva }
\medskip
{\footnotesize
 %% please put the address of the second author
 \centerline{CIMA-UE and Scientific Area of Mathematics, Instituto Superior de Engenharia de Lisboa}
   \centerline{Rua Rom\~{a}o Ramalho, 59, 7000-671 \'{E}vora, Portugal, email lfs@dec.isel.ipl.pt}
} %

\bigskip

%% The name of the associate editor will be entered by an editorial staff
% \centerline{(Communicated by the associate editor name)}

\begin{abstract}
We describe the Williams zeta functions and the twist zeta functions of sub-Lorenz templates generated by renormalizable Lorenz maps, in terms of the corresponding zeta-functions of the sub-Lorenz templates generated by the renormalized map and by the  map that determines the renormalization
type.

\end{abstract}

\section{Introduction}

Let $\phi _t$ be a flow on $S^3$ with
countably many periodic orbits $(\tau _n)_{n=1}^{\infty}$. We can
look to each closed orbit as a knot in $S^3$. It was R. F.
Williams, in 1976, who first conjectured that non trivial knotting
occur in the Lorenz system (\cite{Wi1}). In 1983, Birmann and
Williams introduced the notion of template, in order to study the
knots and links (i.e. finite collections of knots, taking into
account the knotting between them) contained in the geometric
Lorenz attractor (\cite{BW}).

\begin{figure}[ht]\label{charts}\center
  % Requires \usepackage{graphicx}
  \epsfig{file=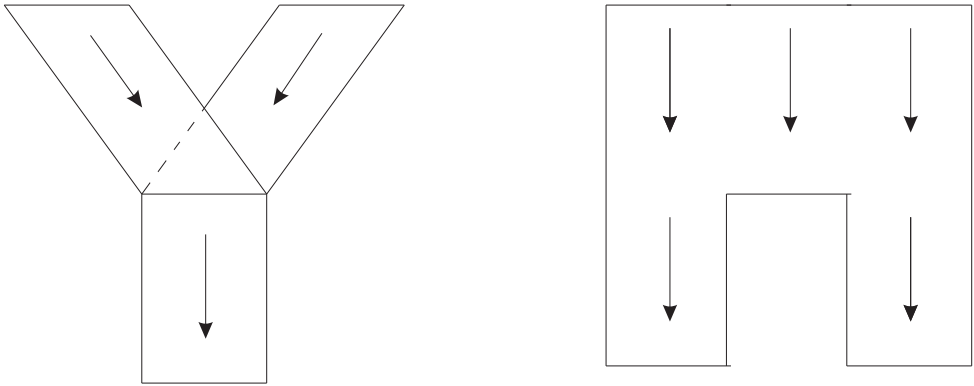,height=2in}\\
  \caption{Charts of templates: joining (left) and splitting (right)}
\end{figure}

A template, or knot holder, consists of a branched two manifold
with charts of two specific types, joining and splitting, together
with an expanding semiflow defined on it, see Figure \ref{charts}.
The relationship between templates and links of closed orbits in
three dimensional flows is expressed in the following result,
known as Template Theorem, due to Birman and Williams in
\cite{BW}.

\begin{theorem}Given a flow $\phi _t$ on a
three-manifold $M$, having a hyperbolic chain-recurrent set, the
link of closed orbits $L_{\phi}$ is in bijective correspondence
with the link of closed orbits $L_{\mathcal{T}}$ on a particular
embedded template $\mathcal{T} \subset M$. On any finite sublink, this
correspondence is via ambient isotopy.
\end{theorem}

We define a \textit{Lorenz flow} as a semi-flow that has a
singularity of saddle type with a one-dimensional unstable
manifold and an infinite set of hyperbolic periodic orbits, whose
closure contains the saddle point (see \cite{MeMA}). A Lorenz
flow, together with an extra geometric assumption (see \cite{Wi})
is called a \textit{geometric Lorenz flow}. The dynamics of this
type of flows can be described by  the iteration of one-dimensional first-return
maps $f:[a,b]\setminus \{c\}\rightarrow [a,b]$ with one discontinuity at $c\in ]a,b[$, increasing in the continuity intervals $[a,c[$ and $]c,b]$ and boundary anchored (i.e. $f(a)=a$ and $f(b)=b$), see \cite{MeMA}. These maps are called Lorenz maps and sometimes we denote them by $f=(f_-,f_+)$, where $f_-$ and $f_+$ correspond, respectively, to the left and right branches.

In \cite{H}, Holmes studied families of iterated horseshoe knots
which arise naturally associated to sequences of period-doubling
bifurcations of unimodal maps.

 It is well known, see for example \cite{DV}, that
period doubling bifurcations in the unimodal family are directly
related with the creation of a 2-renormalization interval, i.e. a
subinterval $J\subset I$ containing the critical point, such that
$f^2|_J$ is unimodal.

Basically there are two types of bifurcations in Lorenz maps (see
\cite{Pi}): one is the usual saddle-node or tangent bifurcations, when
the graph of $f^n$ is tangent to the diagonal $y=x$, and one
attractive and one repulsive $n$-periodic orbits are created or
destroyed; the others are homoclinic bifurcations, when
$f^n(0^{\pm})=f^{n-1}(f_{\pm}(0))=0$ and one attractive
$n$-periodic orbit is created or destroyed in this way, these
bifurcations are directly related with homoclinic bifurcations  of
flows modelled by this kind of maps (see \cite{Pi}).

Considering a  monotone family of Lorenz maps (see \cite{MeMA}), the homoclinic
bifurcations are realized in some lines in the parameters space,
called \textit{hom-lines} or \textit{bifurcation bones}.

It is known that (see \cite{Pi} and \cite{LSSR}), in the context
of Lorenz maps, renormalization intervals are created in each
intersection of two hom-lines. These points are called
\textit{homoclinic points} and are responsible for the
self-similar structure of the bifurcation skeleton of monotone
families of Lorenz maps. So it is reasonable to say that
\textit{homoclinic points are the
 Lorenz version of period-doubling bifurcation points}.

In \cite{Wi2}, Williams considered Lorenz maps with a double saddle connection, i.e., such that $f^n(0^{\pm})=f^{n-1}(f_{\pm}(0))=0$,  and introduced a new determinant-like invariant to classify templates generated by them, in \cite{GHS} these templates are called \textit{sub-Lorenz templates} (note that, in a monotone family of Lorenz maps, maps corresponding to homoclinic points have a double saddle connection). This invariant is a zeta-function that counts periodic orbits taking into account their knot type and the linking between them. In \cite{Su}, Michael Sullivan introduced one other zeta function that counts periodic orbits in positive templates, taking into account the number of twists in their stable manifolds.

In this paper we will study the effect of renormalization over these two invariants, in the context of sub-Lorenz templates. 
\section{Symbolic dynamics of Lorenz maps}

Symbolic dynamics is a very useful combinatoric tool to study the dynamics of one-dimensional maps.

Let $f^j=f\circ f^{j-1}$, $f^0=id$, be the $j$-th iterate of the
map $f$. We define the \textit{itinerary} of a point $x$ under a
Lorenz map $f$ as $i_f(x)=(i_{f}(x))_j, j=0,1,\ldots$, where
\begin{equation*}
(i_{f}(x))_j=\left\{
\begin{matrix}
L & \text{if} & f^j(x)<0 \\
0 & \text{if} & f^j(x)=0 \\
R & \text{if} & f^j(x)>0
\end{matrix}
\right. .
\end{equation*}

It is obvious that the itinerary of a point $x$ will be a finite
sequence in the symbols $L$ and $R$ with $0$ as its last symbol,
if and only if $x$ is a pre-image of $0$ and otherwise it is one
infinite sequence in the symbols $L$ and $R$. So we
consider the symbolic space $\Sigma$ of sequences $X_{0} \cdots
X_{n}$ on the symbols $\{L,0,R\}$, such that $X_{i} \neq 0 $ for
all $i<n$ and: $n=\infty$ or $X_{n}=0$, with the lexicographic
order relation induced by $L<0<R$.

It is straightforward to verify that, for all $x,y \in [-1,1]$, we
have the following:
\begin{enumerate}
\item If $x<y$ then $i_f(x)\leq i_f(y)$, and
\item If $i_f(x)<i_f(y)$ then $x<y$.
\end{enumerate}

We define the \textit{kneading invariant} associated to a Lorenz
map $f=(f_-,f_+)$, as
\begin{equation*}
K_{f}=(K_{f}^{-},K_{f}^{+})=(Li_f(f_-(0)),Ri_f(f_+(0))).
\end{equation*}%

We say that a pair $(X,Y)\in \Sigma \times \Sigma$ is \textit{admissible} if $%
(X,Y)=K_f$ for some Lorenz map $f$. Denote by $\Sigma ^+$, the set of all admissible pairs.

Consider the \textit{shift map} $s:\Sigma\setminus \{0\}
\rightarrow \Sigma$, $s(X_0\cdots X_n)=X_1\cdots X_n$. The set of
admissible pairs is characterized, combinatorially, in the
following way (see \cite{LSSR}).

\begin{proposition}
\label{t3} Let $(X,Y)\in \Sigma \times \Sigma $, then $(X,Y)\in \Sigma ^+$ if and only if $X_{0}=L$, $Y_{0}=R$ and, for $Z\in \{ X,Y \}$ we
have:
\newline (1) If $Z_{i}=L$ then $s ^{i}(Z)\leq X$;
\newline (2) If $Z_{i}=R$ then $s ^{i}(Z)\geq Y$;
with inequality (1) (resp. (2)) strict if $X$ (resp. $Y$) is
finite.
\end{proposition}

On the other hand, the kneading invariant of a map characterizes completely its combinatorics, more precisely, considering a pair $(X,Y)\in \Sigma ^+$, denote by $\Sigma ^+(X,Y)$ the set of sequences $Z\in \Sigma$  that satisfy conditions (1) and (2) from the previous proposition, then we have the following proposition whose proof can be found in \cite{LSSR}.

\begin{proposition}
 Let $X\in \Sigma$ and $f$ be a Lorenz map, then there exists $x\in I$ such that $X=i_f(x)$ if and only if $X\in \Sigma ^+ (K(f))$.
\end{proposition}

\subsection{Renormalization and $*$-product}

In the context of Lorenz maps, we define renormalizability on the
following way, see for example \cite{Pi}:

\begin{definition} Let $f$ be a Lorenz map, then we say
that $f$ is renormalizable if there exist $n,m\in \mathbb{N}$ with $n+m\geq 3$ and  points
$P<y_{L}<0<y_{R}<Q$ such that
$$
g(x)=\left\{
\begin{array}{ll}
f^n(x) & \: if \: y_{L}\leq x<0 \\
f^m(x) & \: if \: 0<x\leq y_{R}%
\end{array}
\right.
$$
is  a Lorenz map.

The map $R_{(n,m)}(f)=g=(f^n,f^m)|_{[y_L,y_R]}$ is called the
$(n,m)$-renormalization of $f$.
\end{definition}

Let $|X|$ be the length of a finite sequence $X=X_0\cdots
X_{|X|-1}0$,  it is reasonable to identify each finite  sequence $X_0\cdots
X_{|X|-1}0$ with the corresponding infinite periodic sequence
$(X_0\cdots X_{|X|-1})^{\infty}$, this is the case, for example,
when we talk about the knot associated to a finite sequence.

Denote $\overline{X}=X_0\ldots X_{|X|-1}$.

It is easy to prove that a pair of finite sequences
$$(X_0\ldots X_{|X|-1}0, Y_0\ldots Y_{|Y|-1}0)$$ is admissible if
and only if the pair of infinite periodic sequences $$((\overline{X})^{\infty},(\overline{Y})^{\infty})$$ is
admissible.

We define the $\ast $-product between a pair of
finite sequences $(X,Y)\in \Sigma \times \Sigma $, and a sequence
$U\in \Sigma $ as
$$
(X,Y)\ast U=\overline{U}_{0}\overline{U}_{1}\cdots
\overline{U}_{|U|-1}0,
$$
where
$$
\overline{U}_{i}=\left\{
\begin{array}{ll}
\overline{X} &\: if \: U_{i}=L \\
\overline{Y} &\: if \: U_{i}=R
\end{array}
\right. .
$$
Now we define the $\ast $-product between two pairs of sequences, $(X,Y),(U,T)\in \Sigma \times \Sigma $, $X$ and $Y$ finite, as
$$
(X,Y)\ast (U,T)=((X,Y)\ast U,(X,Y)\ast T).
$$

The next theorem states  that the reducibility relative to the
$\ast $-product is equivalent to the renormalizability of the map.
The proof can be found, for example, in  \cite{LSSR}.

\begin{theorem}
Let $f$ be a Lorenz map, then $f$ is renormalizable with renormalization $R_{(n,m)}(f)$ iff
there exist two admissible pairs $(X,Y)$ and $(U,T) $ such that
$|X|=n$, $|Y|=m$, $K_f=(X,Y)\ast (U,T)$ and
$K_{R_{(n,m)}(f)}=(U,T)$.
\end{theorem}

We know from \cite{LSSR} that  $(X,Y)\ast (U,T)\in \Sigma ^+$
 if and only if both $(X,Y)\in \Sigma ^+$ and $(U,T)\in \Sigma ^+$, so for each finite admissible pair $(X,Y)$, the
subspace $(X,Y)* \Sigma ^+$ is isomorphic to
the all space $ \Sigma ^+$, this provides a
self-similar structure in the symbolic space of kneading
invariants. It is straightforward to verify that the $*$-product of kneading invariants is associative, consequently this self-similar structure is nested. Now we will state a series of properties concerning the $*$-product, that are in the basis of our results.

The following lemma states that the order structure is reproduced at each level of renormalization.

\begin{lemma}\label{l0}
Let $(X,Y)$ be one admissible pair of finite sequences, and $Z<Z'$,
then $(X,Y)*Z<(X,Y)*Z'$.
\end{lemma}

The proof is straightforward.

For any sequence $X$ and $o\leq p<q < |X|$, we denote 
$$X_{[p,q]}=X_p\ldots X_q.$$

Since we have an order structure in  $\Sigma$, we will denote 
$$[A,B]=\left\lbrace X\in \Sigma : A\leq X \leq B \right\rbrace .$$

For any finite sequence $Z$, consider also the numbers
$$n_L(Z)=\sharp \left\lbrace i\leq |Z|-1 : Z_i=L\right\rbrace $$
and
$$n_R(Z)=\sharp \left\lbrace i\leq |Z|-1 : Z_i=R\right\rbrace .$$
\begin{lemma}\label{l1}
Let $(X,Y)$ and $(S,W)$ be admissible pairs and $A$ and $B$ be any two sequences in $\Sigma $ such that $A\leq B$. Consider $Z\in \left\lbrace  X,Y\right\rbrace $, then a sequence $K\in \Sigma\backslash \left\lbrace Z_{[p,|Z|-1]}0 \right\rbrace $ belongs to $[Z_{[p,|Z|-1]}(X,Y)*A,Z_{[p,|Z|-1]}(X,Y)*B]\cap \Sigma^+ ((X,Y)*(S,W))$ if and only if $K=Z_{[p,|Z|-1]}(X,Y)*C$, with $C\in [A,B]\cap \Sigma^+ ((S,W))$.
\end{lemma}

\textbf{Proof}

Obviously $K_{[0,|Z|-p-1]}=Z_{[p,|Z|-1]}$.

If $K_{|X|-p}=L$, then, since $K\geq Z_{[p,|Z|-1]}(X,Y)*A$, we have that  $A_0=L$ and  $Z_{[p,|Z|-1]}(X,Y)*A=Z_{[p,|Z|-1]}X_{[0,|X|-1]}(X,Y)*\sigma(A)$, so $K_{[|X|-p,2|X]-p-1]} \geq X_{[0,|X|-1]}$.

On the other hand, since $K\in \Sigma ^+((X,Y)*(S,W))$, then $K_{|X|-p}=L$ implies that $K_{[|X|-p,2|X]-p-1]} \leq X_{[0,|X|-1]}$.

Analogously we see that, if $K_{|X|-p}=R$ then $K_{[|X|-p,|X]+|Y|-p-1]} = Y_{[0,|Y|-1]}$.

Repeating these arguments ad infinitum we prove that $K=Z_{[p,|Z|-1]}(X,Y)*W$.

We will now prove that $C\in \Sigma^+((S,W))$.

Let us suppose by contradiction that $C\notin \Sigma^+((S,W))$. Then it happens one of the following situations:
\begin{enumerate}
 \item [(i)] There exists $l$ such that $C_l=L$ and $\sigma ^l(C) >S$.
\item [(ii)] There exists $l$ such that $C_l=R$ and $\sigma ^l(C) <W$.
\end{enumerate}

If it happens situation (i), then there exists $r$ such that $C_{[l,l+r-1]}=S_{[0,r-1]}$ and $C_{p+r}>S_r$. But then, with $q=|X|-p-1+|X|n_L(C_{[0,l-1]})+|Y|n_R(C_{[0,l-1]})$, we have that $K_q=X_0=L$ and $\sigma ^q(K)=\overline{Z}_l\ldots \overline{Z}_{l+r-1}Z_{l+r}\ldots > \overline{S}_0\ldots \overline{S}_{r-1}S_{l+r}\ldots=(X,Y)*S$ and this implies that $K\notin \Sigma ^+((X,Y)*(S,W))$.

If it happens situation (ii), we obtain the contradiction analogously

$\square$

\begin{lemma}\label{l10}
Let $(X,Y)$ be one admissible pair of finite sequences, $0<q<|Y|$
and $Y_q=R$, then $Y_{[q,|Y|-1]}(X,Y)*Z \geq (\overline{Y})^{\infty}$,
for any sequence $Z$. Analogously, if $0<q<|X|$ and $X_q=L$, then
$X_{[q,|X|-1]}(X,Y)*Z \leq (\overline{X})^{\infty}$, for any sequence
$Z$.
\end{lemma}
\textbf{Proof}

Since $(X,Y)$ is admissible, then $Y_{[q,|Y|-1]}(\overline{Y})^{\infty}>(\overline{Y})^{\infty}$, so there exists $l$ such that
$Y_{[q,q+l-1]}=Y_{[0,l-1]}$ and $Y_{q+l}>Y_{l}$ . If
$q+l<|Y|$ the result follows immediately. If $q+l\geq|Y|$, then
necessarily $Y_{|Y|-q}=L$, because otherwise we would have
$(\overline{Y})^{\infty}>Y_{|Y|-q}\ldots$ and $Y_{|Y|-q}=R$, and this violates
admissibility. But then,
$$
Y_{[|Y|-q,|Y|-1]}(\overline{Y})^{\infty}\leq
(\overline{X}) \leq (X,Y)*Z
$$
and this gives the result. The proof of
the second part is analogous. 

$\square$

\begin{lemma}\label{l2}
Let $(X,Y)$ be one admissible pair of finite sequences and
$W,W'\in \{ X,Y \}$. If $s ^p((\overline{W})^{\infty})< s ^q ((\overline{W'})^{\infty})$ and
$W_{[p,|W|-1]} \neq W'_{[q,|W'|-1]}$ then
$$W_{[p,|W|-1]}(X,Y)*Z \leq W'_{[q,|W'|-1]}(X,Y)*Z'$$ for any sequences $Z,Z'$.
\end{lemma}

\textbf{Proof} The proof is divided in four cases:  $W=X$ and
$W'=Y$; $W=Y$ and $W'=X$; $W=W'=X$ and $W=W'=Y$. We will only
demonstrate specifically the first case, since the others follow
with analogous arguments..

Following the hypotheses, there exists $l$ such that $X_{[p,(p+l-1)\mod |X|]}=Y_{[q,(q+l-1)\mod |Y|]}$ and $X_{(p+l)\mod |X|}<Y_{(q+l)\mod |Y|}$. If
$l<\min \{|X|-p,|Y|-q \}$, then the result follows immediately.

If $|X|-p\leq |Y|-q$ and $X_{[p,|X|-1]}=Y_{[q,q+|X|-p-1]}$, then $Y_{q+|X|-p}=R$, because otherwise we would
have $Y_{q+|X|-p}=L$ and $Y_{[q+|X|-p,|Y|-1]}Y^{\infty} >
X^{\infty}$, and this violates  admissibility of $(X,Y)$. So
$Y_{q+|X|-p}=R$ and, from Lemmas \ref{l0} and \ref{l10},
\begin{equation}\label{eq2}
(X,Y)*Z\leq Y^{\infty} \leq
Y_{[q+|X|-p,|Y|-1]}(X,Y)*Z',
\end{equation} and the result follows.

If $|X|-p\geq |Y|-q$ and $X_{[p,p+|Y|-q-1]}=Y_{[q,|Y|-1]}$, then $X_{p+|Y|-q}=L$, because otherwise we would have
$X_{p+|Y|-q}=R$ and $X_{p+|Y|-q}\cdots < Y^{\infty}$, which
contradicts admissibility of $(X,Y)$. So $X_{p+|Y|-q}=L$ and
\begin{equation}\label{eq3}
X_{[p+|Y|-q,|X|-1]}(X,Y)*Z\leq X^{\infty} \leq
(X,Y)*Z'\end{equation} and the result follows.

$\square$

Let us now introduce some more notations:

For $l\leq p$,
$$Z_{[l,p]}=Z_l\cdots Z_p.$$

$$m(A,B)=\min \left\lbrace k\geq 0 : A_{|A|-1-k}\neq B_{|B|-k-1}\right\rbrace .$$
$$\Sigma (A,B) =\left\lbrace \sigma ^n(A), \sigma ^m(B) :0\leq n< |A|, 0\leq m < |B|\right\rbrace ,$$
and $$\phi_{(A,B)}:\Sigma (A,B) \rightarrow \left\lbrace 1,\ldots ,|A|+|B|\right\rbrace ,$$
is the map that associates to each $X\in \Sigma (A,B)$, the position occupied by $X$ in the lexicographic ordenation of $\Sigma(A,B)$.

For each $1\leq k \leq |S|+|W|$, denote 
$$I_k =\left\lbrace \begin{array}{l}
 \left[ (X,Y)*\phi _{(S,W)}^{-1}(k), (X,Y)*\phi _{(S,W)}^{-1}(k+1)\right] \text{ if }m(X,Y)=0 \\
\\
\left[ X_{[|X|-m(X,Y),|X|-1]}(X,Y)*\phi _{(S,W)}^{-1}(k), X_{[|X|-m(X,Y),|X|-1]}(X,Y)*\phi _{(S,W)}^{-1}(k+1)\right] \\
\text{ if }m(X,Y)\neq 0
\end{array}  \right. 
$$

\begin{remark}\label{main}
From the previous three Lemmas  we can take the following conclusions:
\begin{enumerate}

\item If $p<m(X,Y)$, denote by $I_{X_p}$ the set $\left\lbrace X_{[p,|X|-1]}(X,Y)*\sigma ^k(Z) : Z\in \left\lbrace S,W \right\rbrace \: and \: Z_{k-1}=L\right\rbrace $, from Lemmas \ref{l0} and \ref{l1}, $I_{X_p}=\left[ X_{[p,|X|-1]}(X,Y)*\sigma ^2(W),X_{[p,|X|-1]}(X,Y)*\sigma (S)\right] \cap \Sigma((X,Y)*(S,W))$,
analogously denoting  $I_{Y_p}=\left\lbrace Y_{[p,|Y|-1]}(X,Y)*\sigma ^k(Z) : Z\in \left\lbrace S,W \right\rbrace \: and \: Z_{k-1}=R\right\rbrace $ the setwe have that $I_{Y_p}=\left[ Y_{[p,|Y|-1]}(X,Y)*\sigma (W),Y_{[p,|Y|-1]}(X,Y)*\sigma ^2(S)\right] \cap \Sigma((X,Y)*(S,W))$
On the other hand, if $p\geq m(X,Y)$, then $X_{[p,|X|-1]}=Y_{[p,|Y|-1]}$ and \newline
 $\left[ X_{[p,|X|-1]}(X,Y)*\sigma (W) , X_{[p,|X|-1]}(X,Y)*\sigma (S)\right] \cap \Sigma ((X,Y)*(S,W))=\left\lbrace X_{[p,|X|-1]}(X,Y)\sigma ^k(Z) : Z\in \left\lbrace S,W \right\rbrace\right\rbrace $. Without risk of confusion, we will denote these sets by $I_{X_p}$.

\item The ordenation of the elements of the sets $I_{X_p}$ and $I_{Y_q}$ is induced by the ordenation of the sequences $\sigma ^k(Z)$ such that $Z\in \left\lbrace  S,W \right\rbrace $. This follows immediately from Lemma \ref{l0}.

 \item For each $Z\in \{X,Y\}$, if $p\neq |Z|-m(X,Y)-1$ then $\sigma(I_{Z_p})=I_{Z_{p+1}}$. On the other hand, $\sigma(I_{X_{|X|-m(X,Y)-1}}) \cup \sigma(I_{Y_{|Y|-m(X,Y)-1}})=I_{X_{|X|-m(X,Y)}}$. This follows immediately from the definitions.

\item Let $J_{k}=\left[ \max I_{\phi_{(X,Y)}^{-1}(k)},\min I_{\phi_{(X,Y)}^{-1}(k+1)}\right] $ and $H_k=\left[ \phi_{(X,Y)}^{-1}(k),\phi_{(X,Y)}^{-1}(k+1)\right]$, it follows from Lemma \ref{l2} that
 $\sigma \left(J_{p_{k'}}\right)\subset J_{p_{k}}$ iff $\sigma\left(H_{k'}\right)\subset H_k$. Moreover, from (3) and Lemma \ref{l0}, if $\phi_{(X,Y)}^{-1}(k)\notin \left\lbrace \sigma^{|X|-m(X,Y)-1}(X),\sigma^{|Y|-m(X,Y)-1}(Y)\right\rbrace $,
then $\sigma\left( \max I_{\phi _{(X,Y)}^{-1}(k)}\right) = \max I_{\sigma\left( \phi _{(X,Y)}^{-1}(k)\right)} $ 
and $\sigma\left( \min I_{\phi _{(X,Y)}^{-1}(k)}\right) = \min I_{\sigma\left( \phi _{(X,Y)}^{-1}(k)\right)} $
\item From (2) and (3), performing some straightforward computations with the lengths of $(X,Y)*S$ and $(X,Y)*W$ we see that, for $k\neq k'$, then $\sigma ^n(I_k)\cap \sigma ^m(I_{k'})\neq \emptyset$ if and only if both $\sigma ^n(I_k)$ and $\sigma ^m(I_{k'})$ are contained in $I_{|X|-m(X,Y)}$ and $\sigma ^p\left( P_k \right) \cap \sigma ^q\left( P_{k'} \right) \neq \emptyset$, where $P_k=\left[ \phi_{(S,W)}^{-1}(k), \phi_{(S,W)}^{-1}(k+1)\right]$ and $p$ and $q$ are such that $n=|X|n_{L[0,p-1]}\left( \phi_{(S,W)}^{-1}(k)\right) + |Y|n_{R[0,p-1]}\left( \phi_{(S,W)}^{-1}(k)\right)$ and $m=|X|n_{L[0,q-1]}\left( \phi_{(S,W)}^{-1}(k')\right) + |Y|n_{R[0,q-1]}\left( \phi_{(S,W)}^{-1}(k')\right)$.
\end{enumerate}
\end{remark}

This remark is the kernel of the proofs of all our main results, so we will illustrate it with one example.

\begin{example}
 We will consider $(X,Y)=(LRRRL0,RLLR0)$ and $(S,W)=(LRR0,RL0)$. Denoting by $\underset{i}{Z}=\sigma ^i(Z)$, for $Z\in\{X,Y,S,W\}$ and setting that $A<B$ if $A$ is located at the left of $B$, we have the following ordenations of the members of the pairs $(X,Y)$ and $(S,W)$, moreover, the relative position of each member gives the maps $\phi_{(X,Y)}$ and $\phi_{(S,W)}$.
$$\underset{1}{Y} \: \underset{4}{X}\: \underset{2}{Y} \: \underset{0}{X} \: \underset{0}{Y} \: \underset{3}{X} \: \underset{3}{Y} \: \underset{2}{X} \: \underset{1}{X} $$
$$\underset{1}{W} \: \underset{0}{S} \: \underset{0}{W}  \: \underset{2}{S} \: \underset{1}{S} $$  
We will now consider the $*$-product
$$(X,Y)*(S,W)=(\overset{S_0}{\overline{\underset{0}{L}\underset{1}{R}\underset{2}{R}\underset{3}{R}\underset{4}{L}}}\: \overset{S_1}{\overline{\underset{5}{R}\underset{6}{L}\underset{7}{L}\underset{8}{R}}}\: \overset{S_2}{\overline{\underset{9}{R}\underset{10}{L}\underset{11}{L}\underset{12}{R}}}0,\overset{W_0}{\overline{\underset{0}{R}\underset{1}{L}\underset{2}{L}\underset{3}{R}}}\:\overset{W_1}{\overline{\underset{4}{L}\underset{5}{R}\underset{6}{R}\underset{7}{R}\underset{8}{L}}}0),$$
we used the underscripts to indicate the corresponding iterate of the shift map and the upperscripts to indicate the element of $(S,W)$ that generated each subword in the $*$-product. Denoting by $(A,B)=(X,Y)*(S,W)$ we will now order the elements of the $*$-product pair:
$$\overset{I_{Y_1}}{\overleftrightarrow{\underset{1}{B}\:\underset{10}{A}\:\underset{6}{A}}} 
\overset{J_1}{\overbrace{\:\:\:}} \overset{I_{X_4}}{\overleftrightarrow{\underset{8}{B}\:\underset{4}{A}}}\overset{J_2}{\overbrace{\:\:\:}} \overset{I_{Y_2}}{\overleftrightarrow{\underset{2}{B}\:\underset{11}{A}\:\underset{7}{A}}}
\overset{J_3}{\overbrace{\:\:\:}} \overset{I_{X_0}}{\overleftrightarrow{\underset{4}{B}\:\underset{0}{A}}}\: \overset{I_{Y_0}}{\overleftrightarrow{\underset{0}{B}\:\underset{9}{A}\:\underset{5}{A}}}
\overset{J_4}{\overbrace{\:\:\:}} \overset{I_{X_3}}{\overleftrightarrow{\underset{7}{B}\:\underset{3}{A}}}\overset{J_5}{\overbrace{\:\:\:}} \overset{I_{Y_3}}{\overleftrightarrow{\underset{3}{B}\:\underset{12}{A}\:\underset{8}{A}}}
\overset{J_6}{\overbrace{\:\:\:}} \overset{I_{X_2}}{\overleftrightarrow{\underset{6}{B}\:\underset{2}{A}}}\overset{J_7}{\overbrace{\:\:\:}} \overset{I_{X_1}}{\overleftrightarrow{\underset{5}{B}\:\underset{1}{A}}} $$
as we can see, the ordered disposition of the members of $(X,Y)*(S,W)$ is obtained from the ordered disposition of $(X,Y)$, inflating each $\underset{i}{Z}$, $Z\in \{X,Y \}$ and substituting it by the elements of the corresponding $I_{Z_i}$, ordered according with the ordenation of the $\underset{k}{H}$, $H\in \{S,W\}$.
\end{example}

\section{Sublorenz Templates}

Now we will follow \cite{Wi2} to introduce Sub-Lorenz templates and Williams zeta-functions.

 We say that a Lorenz map $f$ haves a double saddle connection if
 $f^n(0^-)=f^m(0^+)=0$ for some $n,m$.

 In this case the points $\{f^i(0^-),f^j(0^+): 1\leq i \leq n,1\leq j \leq m\}$ we can define a finite Markov partition for the
 semiflow.

 The restriction of the semiflow to this partition is called a
 Sub-Lorenz template.

\begin{example}
For $(K^-,K^+)=(LRR0,RL0)$, we construct the sub-Lorenz template $T_{(K^-,K^+}$ following the procedure of Figure \ref{figtempl}
\begin{figure}\label{figtempl}
\begin{center}
\epsfig{file=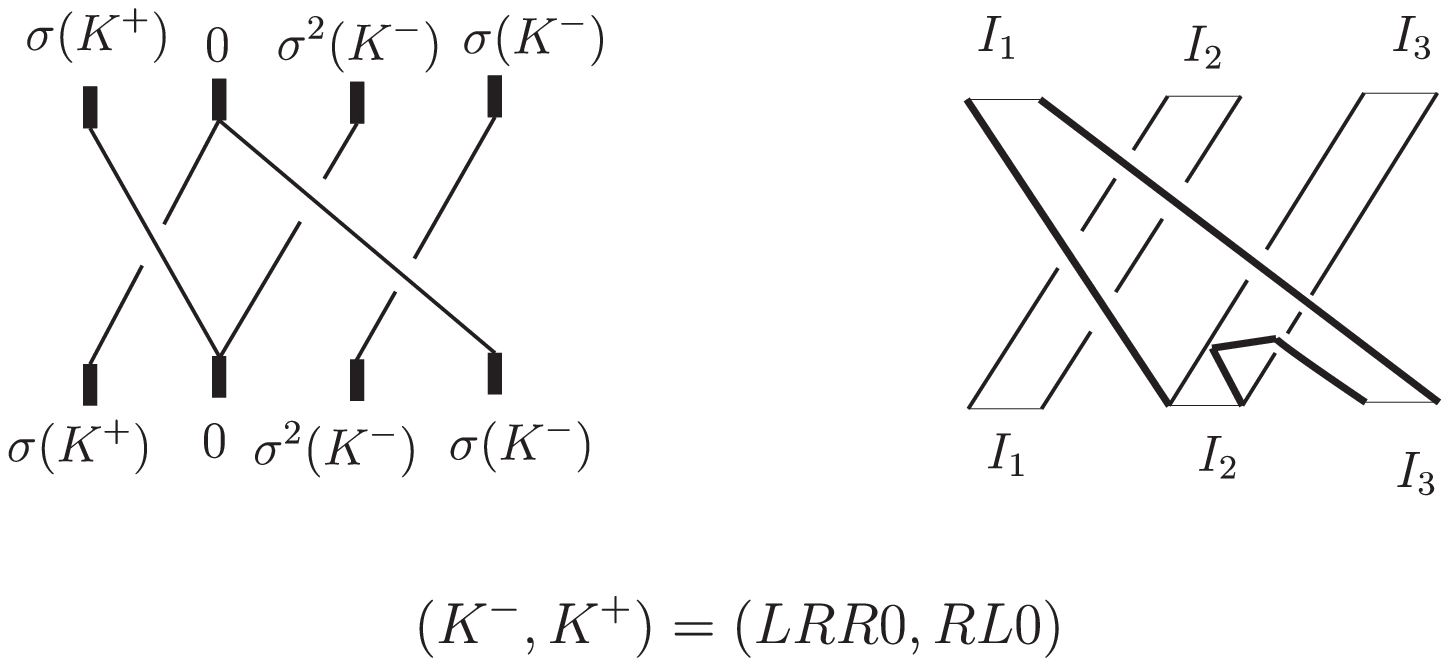,height=1.8in}
\caption{Sub-Lorenz template $T_{(LRR0,RL0)}$.}
\end{center}
\end{figure}

\end{example}

We identify each Sub-Lorenz template with the corresponding kneading invariant $(X,Y)$ and denote it with $T_{(X,Y)}$. We associate to the template the  transition
 matrix $A_{T_{(X,Y)}}=[a_{ij}]$ where
 $$a_{ij}=\left\{\begin{array}{ll}
   1 & \text{ if } I_j\subset f(I_i) \\
   0 & \text{ if } I_j\bigcap f(I_i) =\emptyset \
 \end{array}\right.
$$

Now we associate to $T_{(X,Y)}$ a labeled transition matrix $A_{T_{(X,Y}}(L,R)=[a'_{ij}]$, where
$$a'_{ij}=\left\{\begin{array}{ll}
   L & \text{ if } I_j\subset f(I_i) \text{ and } I_i \text{ is located on the left of }0\\
    R & \text{ if } I_j\subset f(I_i) \text{ and } I_i \text{ is located  on the right of }0 \\
   0 & \text{ if } I_j\bigcap f(I_i) =\emptyset \
 \end{array}\right.
$$

 In the previous example we have
$$A_{T_{(X,Y)}}=\left[\begin{array}{ccc}
  0 & 1 & 1 \\
  1 & 0 & 0 \\
  0 & 1 & 0
\end{array}
\right]
$$
and
$$A_{T_{(X,Y)}}(L,R)=\left[\begin{array}{ccc}
  0 & L & L \\
  R & 0 & 0 \\
  0 & R & 0
\end{array}
\right]
$$
 Given a group $F$ and a matrix $A$ with entries in
$F$, we define the \textit{link-ring} $R(A)$ as follows:
\begin{enumerate}
\item For each sequence (called a cycle below)
$$i_1,i_2,\ldots , i_k, \:\: i_r\neq i_s \text{ for } r\neq s,$$
such that the product
$$a_{i_1i_2}a_{i_2i_3}\cdots a_{i_ki_1}\neq 0,$$
let $(a_{i_1i_2}a_{i_2i_3}\cdots a_{i_ki_1})$ be the equivalence
class under cyclic permutations of this product. These equivalence
classes are called \textbf{\textit{free-knot symbols}} and the
indices $i_1,i_2,\ldots , i_k$ are called \textit{nodes} of the
free knot symbol.
\item Define $R(A)$ to be the free abelian group generated by the
free knot symbols defined in (1).
\end{enumerate}
The nodes in the matrix correspond to the cells in
the Markov partition. The cycles in the matrix correspond to
periodic orbits on the Markov partition, and thus to periodic
orbits in the flow. Since there is no natural way to choose a
specific point for a periodic orbit, the natural invariant is the
cyclic permutation class.

 \textbf{We do not permit products of letters to
commute}, since this usually corresponds to different orbits and
consequently to different knot types, for example $L^3R^2$ and
$L^2RLR$ correspond, in the first case to the unknot and in the
second to the trefoil knot.

However, when we consider a product of words, it
represents the union of two knots (i.e.  \textbf{link}) and in
this case we don't care which happens first, so \textbf{we permit
products of free-knot symbols (words) to commute}.
By a \textbf{\textit{free-link symbol}} in the ring
$R(A)$is meant a product $x_1\cdots x_l$ of free-knot symbols, no
two of which have a node in common.

(Just as a product of two cycles can occur in a determinant, only
if they have no node in common)

In \cite{Wi2}, Williams defined the following determinant-like
function and proved the following two results about it:
\begin{definition}
Let $A_{T_{(X,Y}}(L,R)$ be a labeled transition matrix of a sub-Lorenz template, then 
$$\text{link-det}(I-A_{T_{(X,Y}}(L,R))=\underset{free-link\:
symbols}{\sum}(-1)^lx_1x_2\cdots x_l .$$ 
\end{definition}
\begin{theorem}[Williams]For Lorenz attractors $T_{(X,Y}$ with a double saddle
connection,
$$\text{link-det}(I-A_{T_{(X,Y}}(L,R))=\underset{\cal{L}}{\sum}(-1)^{|L|}fls(L),$$
where $\cal{L}$ is the collection of all links $L$ in the
attractor which have at most one point in each partition set,
$|L|$ is the number of components in $L$, and $fls(L)$ means the
free-link symbol of $L$.
\end{theorem}

\begin{theorem}[Williams] For Lorenz attractors $T_{(X,Y}$ with
a double saddle connection, we
have that
$$exp\left(-\underset{i=0}{\overset{\infty}{\sum}}tr(A_{T_{(X,Y}}(L,R)^i)/i\right)=\text{link-det}(I-A_{T_{(X,Y}}(L,R)).$$
\end{theorem}
From this result, we name $\zeta _W(T_{(X,Y})=\text{link-det}(I-A_{T_{(X,Y}}(L,R))$, the
\textbf{\textit{Williams zeta function}} of the template $T_{(X,Y}$.

Finally we state our factorization result about Williams zeta-functions:

\begin{theorem} Let $T_{(X,Y)*(S,W)}$ be a Sub-Lorenz template generated by a Renormalizable Lorenz map $f$, with a double saddle connection and  kneading invariant
kneading invariant $K(f)=(X,Y)*(S,W)$, then
$$\zeta_W(T_{(X,Y)*(S,W)})+1=\left[\zeta_W(T_{(X,Y)}) +1 \right] \times \left[(X,Y)*\zeta_W(T_{(S,W)}) +1\right],$$
where $(X,Y)*\zeta_W(T_{(S,W)})= \sum (-1)^l(X,Y)*x_1 \ldots (X,Y)*x_l$ such that the sum is taken over all  free-link symbols of $T_{(S,W)}$, $x_1\ldots x_n$.
\end{theorem}

\textbf{Proof}
Denote $\mathcal{RB}=\underset{Z\in\{X,Y\}}{\cup} \underset{0\leq i<|Z|}{\cup}I_{Z_i}$.

\begin{enumerate}

\item[i] Let $\mathcal{P}_{(X,Y)*(S,W)}$ be the Markov partition associated to $(X,Y)*(S,W)$, and  $I$ be a cell of $\mathcal{P}_{(X,Y)*(S,W)}$, then $I\subset I_Z$ for some $Z\in  \mathcal{RB}$ or $I=J_{k}$ for some $k$ such that $1\leq k \leq |X|+|Y|$. Moreover, from 3 of remark \ref{main} we have that $\sigma ^n( \mathcal{RB} )\subset \mathcal{RB}$ for all $n$ and, from 4, the free-knot symbols associated to nodes in $\mathcal{P}_{(X,Y)*(S,W)}\backslash  \mathcal{RB}$ are exactly the same free-knot symbols from $T_{(X,Y)}$ and consequently the free-link symbols are also the same.

\item[ii] From (5) of Remark \ref{main}, the free-link symbols associated to knots in $\Sigma\mathcal{RB}$ are exactly those that can be written in the form $(X,Y)*x_1 \cdots (X,Y)*x_l$, where $x_1 \cdots x_l$ are free-link symbols of $T_{(S,W)}$.

\item[iii] Finally, from (3) of Remark \ref{main} there are no cycles, simultaneously with nodes in $\mathcal{RB}$ and in $\mathcal{P}_{(X,Y)*(S,W)}\backslash \mathcal{RB}$, so all free-link symbols of type $\alpha \beta$ where $\alpha$ is a free link symbol associated to nodes in $\mathcal{RB}$ and $\beta$ is a free link symbol associated to nodes in $\mathcal{P}\backslash \mathcal{RB}$ will be terms of the zeta-function of $T_{(X,Y)*(S,W)}$.

\end{enumerate}

Finally, from (i), (ii) and (iii) we have that 
$$\zeta _W(T_{(X,Y)*(S,W)})=\zeta _W(T_{(X,Y)})+(X,Y)*\zeta _W(T_{(S,W)})+(X,Y)\times (X,Y)*\zeta _W(T_{(S,W)}).$$

$\square$

\section{Twist-zeta function}
In this section we will present a factorization formula for the Twist-zeta function presented by Sullivan in \cite{Su}.

A ribbon is an embedded annulus or M\"obius band in $S^3$. Like knots and templates, ribbons can be braided. A ribbon which has a braid presentation such that each crossing of one strand  over another is positive and each twist in each strand is positive, will be called a positive ribbon. The core and boundary of positive ribbons are positive braids. 

\begin{definition}
If   $R$ is a ribbon and $b(R)$ is a braid presentation of $R$, we  define the computed twist
$$\tau _c(R)=2n+t,$$
where $t$ is the sum of the half twists in the strands of $b(R)$ and $n$ is the number of strands of the core.
\end{definition}

In \cite{Su}, Sullivan proved that $\tau _c$ is an isotopy invariant of positive ribbons over positive braid presentations, so the definition is consistent. 

\begin{definition}
Given a template, $T$, and  an orbit, $O$, on $T$, we define the ribbon  $R(T,O)$,  as the ribbon defined by the unit normal bundle of $O$.  
\end{definition}

Sullivan proved that, for positive templates, the number of closed orbits with a given computed twist  is finite. This permited him to formulate the following definition:

\begin{definition}
 For a given template $T$, let $T_{q'}$ be the number of closed orbits with computed twist $q'$. Let $\mathcal{T}_q=\sum_{q'|q}q'T_{q'}$. Define the Sullivan zeta function of the template to be the exponential of a formal power series:
$$\zeta^S _T (t)=\exp \left(\underset{q=2}{\overset{\infty}{\sum}}\mathcal{T}_q\dfrac{t^q}{q}\right) .$$
\end{definition}

We now define a twist matrix $A(t)=[a_{ij}]$ whose entries are nonnegative powers of $t$ and $0$'s, by considering the contribution to $\tau _c$ as an orbit goes from one element of a Markov partition to another. Let $a_{ij}=0$ if there is no branch going from the $i$-th to the $j$-th partition element and $a_{ij} = t^{q_{ij} }$ if there is such a branch, where 
$q_{ij}$ is the amount of computed twist an orbit picks up as it travels from the $i$-th to the $j$-th partition element.

In the case of a sub-Lorenz template $T_{(X,Y)}$, since we have one curl in each ribbon from one partition element to other, see Figure \ref{figtempl} then we obtain $A_{T(X,Y)}(t)$ simply substituting each $1$ element of the transition matrix by $t^2$, i.e.
$$A_{T(X,Y)}(t)=A_{T(X,Y)}(t^2,t^2).$$

Sullivan proved the following theorem in \cite{Su}.

\begin{theorem}
For any template $T$ and any allowed choice of  $A(t)$ we have 
$$\zeta ^S_T(t)=\dfrac{1}{det(I-A(t))}.$$
\end{theorem}

For a sub-Lorenz template $T_{(X,Y)}$, denote 
$$\zeta ^S_{T_{(X,Y)}}(L,R)=\dfrac{1}{det(I-A_{T_{(X,Y)}}(L,R))}.$$
Following this notation, we have that $\zeta ^S_{T_{(X,Y)}}(t)=\zeta ^S_{T_{(X,Y)}}(t^2,t^2)$.

We will now state our factorization result.

\begin{theorem}
 For a reducible kneading pair $(X,Y)*(S,W)$ with both $(X,Y)$ and $(S,W)$, admissible finite Lorenz pairs we have that
$$\zeta ^S_{T_{(X,Y)*(S,W)}}(t^2,t^2)=\zeta ^S_{T_{(X,Y)}}(t^2,t^2)\zeta ^S_{T_{(S,W)}}(t^{2|X|},t^{2|Y|}).$$

\end{theorem}

\textbf{Proof}
From Remark \ref{main}, the Markov partition $\mathcal{P}$ associated to $T_{(X,Y)*(S,W)}$ can be splited as $\mathcal{P} =\mathcal{RB}\cup \mathcal{P} \backslash \mathcal{RB}$ in such a way that all the iterates of each periodic orbit are exclusively contained in $\mathcal{RB}$ or in $\mathcal{P} \backslash \mathcal{RB}$. So we can split the sum
$$\underset{q=2}{\overset{\infty}{\sum}}\mathcal{T}_q\dfrac{t^q}{q}=\underset{q=2}{\overset{\infty}{\sum}}\mathcal{T}_q\left(\mathcal{RB}\right)\dfrac{t^q}{q}+\underset{q=2}{\overset{\infty}{\sum}}\mathcal{T}_q\left(\mathcal{P}\backslash\mathcal{RB}\right)\dfrac{t^q}{q}$$
where $\mathcal{T}_q\left(\mathcal{RB}\right)$ (respectively $\mathcal{T}_q\left(\mathcal{P}\backslash\mathcal{RB}\right)$) means simply that we are counting orbits contained in $\mathcal{RB}$ (respectively in $\mathcal{P} \backslash \mathcal{RB}$).

From (4) of Remark \ref{main}, $\exp \left(\underset{q=2}{\overset{\infty}{\sum}}\mathcal{T}_q\left(\mathcal{P}\backslash\mathcal{RB}\right)\dfrac{t^q}{q}\right) =\zeta ^S_{T_{(X,Y)}}(t^2,t^2)$. Now, from (5) of Remark \ref{main}, each ribbon leaving a cell $I_k$ makes $|X|$ curls without splitting if $I_k$ is on the left of $0$ and makes $|Y|$ curls without splitting if $I_k$ is on the left of $0$, before reenter in $I_{X_{|X|-m(X,Y)-1}}$. Moreover $\sigma^{|X|}(I_k)\cap I_{k'}$ (respectively $\sigma^{|Y|}(I_k)\cap I_{k'}$) is not empty if and only if $\sigma(P_k) \cap \sigma (P_{k'}) \neq \emptyset $, so $\exp \left( \underset{q=2}{\overset{\infty}{\sum}}\mathcal{T}_q\left(\mathcal{RB}\right)\dfrac{t^q}{q}\right)=\zeta ^S_{T_{(S,W)}}(t^{2|X|},t^{2|Y|})$ and the result follows.

$\square$

\end{document}